\documentclass[10pt]{article}
\usepackage{amsmath, amssymb}
\usepackage{amsthm} 
\usepackage[all]{xy}

\theoremstyle{plain} 

\theoremstyle{definition}

\makeatletter
\def\even{\text{\rm even}}
\def\odd{\text{\rm odd}}

\def\L{\text{{\cal L}}}
\def\TT{TM\oplus T^*M}
\def\GCS{generalized complex structure}
\def\Int{\text{\rm Int}}
\title{Unobstructed deformations of generalized complex structures induced by $C^\infty$ logarithmic symplectic structures \\and logarithmic Poisson structures} 

\author{  Ryushi Goto \\ 
\small Department of Mathematics,
\small Graduate School of Science, Osaka University \\
\small goto@math.sci.osaka-u.ac.jp }

\date{} 
\begin{document}
\maketitle

\footnote{ 
\textit{Key words and phrases}. 
}
\footnote{ 
$^{*}$Partly supported by the Grant-in-Aid for Scientific Research (B),
Japan Society for the Promotion of Science. 
} 
%
%
%

\def\GCS{\text{\rm generalized complex structure}}

\def\Sch{\text{\rm Sch}}
\def\F{\cal F}
\def\L{{\cal L}}
\def\I{{\cal I}}
\def\bgn{\begin}
\def\Lrarrow{\Leftrightarrow}
\def\CL{\text{\rm CL}}
\def\CY{\text{\scriptsize\rm CY}}
\def\GL{\text{\rm GL}}
\def\E{\bold E}
\def\J{{\cal J}}
\def\L{{\cal L}}
\def\Ob{\text{\rm Ob}}
\def\1{{[1]}}
\def\2{{[2]}}
\def\3{{[3]}}
\def\({\left(}
\def\){\right)}
\def\s-circ{\,{\scriptstyle{\circ}}\,}
\def\<<{<\negthinspace \negthinspace<}
\def\Ad{\text{\rm Ad}}
\def\ad{\text{\rm ad}}
\def\noin{\noindent}
\def\even{\text{\rm even}}
\def\bgn{\begin}
\def\aln{\bgn{align}}
\def\endaln{\end{align}}
\def\cou{\ss\text{\rm co}}
\def\bE{\bold E}
\def\<{<\negthinspace \negthinspace <}
\def\ssw{{\scriptscriptstyle\wedge}}
\def\t{\theta}
\def\TX{TX\oplus T^*X}
\def\TM{TM\oplus T^*M}
\def\({\left(}
\def\){\right)}
\def\Im{\text{\rm Im}}
\def\Re{\text{\rm Re}}
\def\[{\big[\neg\big[}
\def\]{\big]\neg\big]}
\def\al{\al}
\def\cpt{\text{\rm cpt}}
\def\M{{\cal M}}
\def\ind{\text{\rm ind}}
\def\diam{\text{\rm diam}}
\def\Vol{\text{\rm Vol}}
\def\tr{\text{\rm tr}}
\def\a{\alpha}
\def\b{\beta}
\def\bs{\backslash}
\def\e{\varepsilon}
\def\gam{\gamma}
\def\Gam{\Gamma}
\def\k{\kappa}
\def\del{\delta}
\def\lam{\lambda}
\def\ome{\omega}
\def\Ome{\Omega}
\def\sig{\sigma}
\def\Sig{\Sigma}
\def\eps{\varepsilon}
\def\thet{\theta}
\def\The{\Theta}
\def\A{{\cal A} }
\def\E{{\cal E}}
\def\f{\bold {f}}
\def\G{G_2}
\def\D{\Bbb D}
\def\Diff{\text{\rm Diff}_0}
\def\Q{\Bbb Q}
\def\R{\Bbb R}
\def\C{\Bbb C}
\def\H{\text{\rm H}}
\def\Z{\Bbb Z}
\def\O{\Bbb O}
\def\P{\frak P}
\def\M{\frak M}
\def\N{\frak N}
\def\V{{\scriptscriptstyle V}}
\def\X{{\cal X}}
\def\w{\wedge}
\def\({\left(}
\def\){\right)}
\def\G{G_2} 
\def\Cl{\frak Cl} 
\def\nab{\nabla}
\def\neg{\negthinspace}
\def\h{\hat}
\def\hrho{\hat{\rho}}
\def\wideh{\widehat}
\def\til{\tilde}
\def\wtil{\widetilde}
\def\ch{\check}
\def\l{\left}
\def\ol{\overline}
\def\pa{\partial}
\def\olpa{\ol{\partial}}
\def\r{\right}
\def\ran{\rangle} 
\def\lan{\langle}
\def\Spin{\text{Spin}(7)}
\def\ss{\scriptscriptstyle}
\def\trian{\triangle}
\def\hook{\hookrightarrow}
\def\hyper{hyperK\"ahler}
\def\Ker{\text{\rm Ker}}
\def\arrow{\longrightarrow}
\def\lrarrow{\Longleftrightarrow}
\def\CP{\dsize{\Bbb C} \bold P}
\def\SL{\ss{\text{\rm SL}}}
\def\Sym{\text{Sym}^3\,\h{\Ome}^1 _X}
\def\tilSym{\widetilde{ \text{Sym}^3}\,\h{\Ome}^1_X }
\def\SymS{\text{Sym}^3\,{\Ome}^1_S }
\def\W{\wedge_{\ome^0}}
\def\sw{{\ss{\w}}}
\def\noin{\noindent}
\def\bsh{\backslash}
\def\reg{\text{\rm reg }}
\def\:{\, :\,}
\def\CL{\text{\rm CL}}
\def\TT{T\oplus T^*}
\def\complex{generalized complex }
\def\K\"ahler{generalized K\"ahler}
\def\vol{\text{\rm vol}}
\def\MC{\text{\rm Maurer-Cartan equation}}
\def\10{\displaystyle L^{10}}
\def\2{\displaystyle L^2}
\def\c0{\displaystyle C^0}
\def\dstyle{\displaystyle}
\def\10{\displaystyle L^{10}}
\def\2{\displaystyle L^2}
\def\del{\delta}
\def\del2{\displaystyle L^2_{0,\delta}}
\def\c0{\displaystyle C^0}
\def\dstyle{\displaystyle}
\def\Lp{\dstyle L^p}
\def\del{\delta}
\def\cl{\text{\rm cl}}
\def\FOW{\text{\tiny\rm FOW}}
\def\K{{\cal K}}
\def\M-A{\text{\rm Monge-Amp\`ere}}
\def\O{{\cal O}}
\def\cyl{\text{\rm cyl}}
\def\cone{\text{\rm cone}}
\def\Ric{\text{\rm Ric}}
\def\M-A{\text{\rm Monge-Amp\`ere}}
\def\[{\big[\,}
\def\]{\,\big]}
\def\End{\text{\rm End\,}}
\def\Hom{\text{\rm Hom}}
\def\SO{\text{\rm SO}}
\def\id{\text{\rm id}}
\def\P{\Bbb P}
\def\even{\text{\rm even}}
\def\odd{\text{\rm odd}}
\def\Int{\text{\rm Int}}
\abstract{We shall introduce the notion of $C^\infty$ logarithmic symplectic structures on a differentiable manifold which
is an analog of the one of logarithmic symplectic structures in the holomorphic category.
We show that the generalized complex structure induced by a $C^\infty$ logarithmic symplectic structure
has unobstructed deformations  
 which are parametrized by an open set of the second de Rham cohomology group of the complement of 
type changing loci if the type changing loci are smooth. 
Complex surfaces with smooth effective anti-canonical divisors admit unobstructed deformations of 
generalized complex structures such as del pezzo surfaces and Hirzebruch surfaces. 
We also give some calculations of Poisson cohomology groups on these surfaces.
 Generalized complex structures $\J_m$ on the connected sum $(2k-1)\C P^2\# (10k-1)\ol {\C P^2}$ 
as in \cite{Cavalcanti_Gualtieri_2006}, \cite{Goto_Hayano} are induced by 
 $C^\infty$ logarithmic symplectic structures modulo the action of $b$-fields and it turns out that 
 {generalized complex structure}s $\J_m$ have
 unobstructed deformations of dimension $12k+2m-3$.}
\tableofcontents
\section{Introduction}
\label{sec:1}
Generalized complex structures are mixed geometric structures building a bridge between complex geometry and real symplectic geometry. 
Both complex structures and  real symplectic structures give rise to {generalized complex structure}s of special classes. However 
 a {generalized complex structure} on a manifold can admit the type changing loci on which the type of the structure can change form the one from a real symplectic structure to the one from a complex structure.\par
A complex surface $S$ with a non trivial holomorphic Poisson structure $\b$ has a generalized complex structure $\J_\b$ with type changing loci at zeros of $\b$.
It is striking that  $(3k-1)\C P^3\# (10k-1)\ol {\C P^2}$ does not admit complex structures and 
real symplectic structures, but  $(3k-1)\C P^3\# (10k-1)\ol {\C P^2}$ has {generalized complex structure}s $\J_m$ with $m$ type changing loci \cite{Cavalcanti_Gualtieri_2006}, \cite{Goto_Hayano}.

The deformation complex of a {generalized complex structure}  is given  by the Lie algebroid complex $(\w^\bullet \ol L_\J, d_L)$. Then  the space of infinitesimal deformations is the second cohomology group $H^2(\w^\bullet \ol L_\J)$ and  the obstruction space is the third cohomology group $H^3(\w^\bullet\ol L_\J)$.
Then the Kuranishi families of {generalized complex structure}s are constructed \cite{Gualtieri_2011}.\par

{\it A $C^\infty$ logarithmic symplectic structure} $\ome_\C$ on a manifold $M$ along a submanifold $D$ of
 real codimension $2$ is a complex $2$-form which is
given by the following on a neighborhood of $D$,
$$\ome_\C=\frac{dz_1}{z_1}\w dz_2+dz_3\w dz_4+\cdots+dz_{2m-1}\w dz_{2m},$$
where $(z_1,\cdots, z_{2m})$ are complex coordinates and
$D=\{z_1=0\}$. 
On a neighborhood of the complement $M\bsh D$, $\ome_\C$ is a smooth $2$-form $b+\sqrt{-1}\ome$  
where $b$ is a $d$-closed $2$-form and $\ome$ denotes a real symplectic structure.
(Note that we do not assume that $M$ is a complex manifold.)
The exponential $\phi:=e^{\ome_\C}$ gives rise to a generalized complex structure $\J_\phi$

One of the purposes of the paper is to show that the generalized complex structure $\J_\phi$ induced by a $C^\infty$ logarithmic symplectic structure $\ome_\C$ 
has unobstructed deformations which are given by 
the second cohomology group of the complement of type changing loci if type changing loci are smooth 
(see 
Theorem \ref{th:log.cohomology} and Theorem \ref{th:unobst. log GCS}).
Our unobstructedness theorems may be regarded as an analog of the unobstructedness theorems of 
Calabi-Yau and hyperK\"ahler manifolds. However we do not use the Hodge theory and the $\pa\ol\pa$-lemma
but our method is rather topological.

The generalized complex structure $\J_\b$ on a Poisson surface $S$ and 
$\J_m$ on the $(3k-1)\C P^3\# (10k-1)\ol {\C P^2}$ are induced by  $C^\infty$ logarithmic symplectic structures 
modulo the action of $b$-fields. Thus we can apply our unobstructedness theorems to these generalized complex structures
(see Theorem \ref{th:ubobstructed defs of Jb} and Theorem \ref{unobst. def of (2k-1)CP2 (10k-1)ol CP2}).

The paper is organized as follows. 
In Section 2 we give a short explanation of {generalized complex structure}s and in Section 3 we also provide fundamental notions of deformations of {generalized complex structure}s. It is remarkable that Poisson cohomology groups
of a Poisson structure $\b$ which are the hypercohomology groups of the Poisson complex coincide with the Lie algebroid cohomology groups of the generalized complex structure $\J_\b$ induced by $\b$. 
In Section 4, we discuss the {generalized complex structure}s induced from a logarithmic symplectic structure along a smooth divisor $D$ on a complex manifold $X$ which is the dual of a holomorphic Poisson structure $\b$.
In Section 5, we show that the Poisson cohomology groups of $\b$ are isomorphic to the de Rham cohomology groups of the complement $X\bsh D$. 
In Section 6, we introduce $C^\infty$ logarithmic symplectic structures and show our main theorems. 
In Section 7, we discuss logarithmic deformations. 
Proofs of our main Theorems are given in Section 8.
In Section 9, the unobstructedness theorem can be applied to a complex surface with smooth anti-canonical divisor to obtain unobstructed deformations.
We calculate Poisson cohomology groups and the de Rham cohomology groups of the complements. It is intriguing that there are differences between these two cohomology groups when $D$ has singularities (see Remark \ref{re:cubic curve with one node}).
In Section 9.3, we explain a construction of {generalized complex structure}s $\J_m$ on $(3k-1)\C P^3\# (10k-1)\ol {\C P^2}$ by logarithmic transformations with $m$ type changing loci.
We apply our theorems to {generalized complex structure}s on $(3k-1)\C P^3\# (10k-1)\ol {\C P^2}$ to obtain unobstructed deformations of {generalized complex structure}s (see Theorem \ref{unobst. def of (2k-1)CP2 (10k-1)ol CP2}).

\numberwithin{equation}{section}
\section{Generalized complex structures}
Let $TM$ be the tangent bundle on a differentiable manifold of dimension $2n$ and $T^*M$ the cotangent bundle of $M$.
The symmetric bilinear form $\lan\,,\,\ran$ on 
the direct sum $TM \oplus T^*M$ is defined by 
$\lan v+\xi, u+\eta \ran=\frac12\(\xi(u)+\eta(v)\), $ where 
$u, v\in TM, \xi, \eta\in T^*M.$ 
Then the symmetric bilinear form $\lan\,,\,\ran$ yields the fibre bundle 
SO$(TM\oplus T^*M)$ with fibre the special orthogonal group.
A section of bundle SO$(TM\oplus T^*M)$ is an endomorphism of $TM\oplus T^*M$ preserving $\lan\,,\,\ran$ and its determinant is equal to one.
If a section $\J$ of SO$(TM\oplus T^*M)$ satisfies 
$\J^2=-\id,$ then $J$ is called an almost \complex structure which gives the
decomposition 
$(TM \oplus T^*M)^\C =L_\J \oplus \ol L_\J$ into eigenspaces, where 
$L_\J$ is the eigenspace of eigenvalue $\sqrt{-1}$ and $\ol L_\J$ is the complex conjugate of $L_\J$.
The Courant bracket is defined by 
$$
 [u+\xi, v+\eta]_{\cou}=[u,v]+{\cal L}_u\eta-{\cal L}_v\xi-\frac12(di_u\eta-di_v\xi),
 $$
 where $[u,v]$ denotes the bracket of vector fields $u$ and $v$ and 
 ${\cal L}_v\eta$ and ${\cal L}_u\xi$ are the Lie derivatives and 
 $i_u\eta$ and $i_v\xi$ stand for the interior products. 

If $L_\J$ is closed with respect to the Courant bracket, then $\J$ is a \complex structure, that is, 
$[e_1, e_2]_{\cou}\in L_\J$ for all $e_1, e_2\in L_\J$.

The direct sum $TM\oplus T^*M$ acts on differential forms by the 
 interior and exterior products,
$$
e\cdot\phi=(v+\eta)\cdot\phi =i_v\phi+\eta\w\phi,
$$ 
where $e=v+\eta$ and $v\in TM$ and $\eta\in T^*M$ and $\phi$ is a differential forms.
Then it turns out that 
$$
e\cdot e\cdot\phi=(v+\eta)\cdot(v+\eta)\cdot\phi=\eta(v)\, \phi=\lan e, e\ran\phi.
$$
Since it is the relation of the Clifford algebra with respect to $\lan\,,\,\ran$, 
we obtain the action of the Clifford algebra bundle Cl$(TM\oplus T^*M)$ on differential forms which is the spin representation.

We define $\ker\Phi:=\{ E\in (TM\oplus T^*M)^\C\, |\, E\cdot\Phi=0\, \}$ for a differential form $\Phi
\in \w^{\even/\odd}T^*M.$
If $\ker\Phi$ is maximal isotropic, i.e., $\dim_\C\ker\Phi=2n$, then $\Phi$ is called {\it a pure spinor} of even/odd type.

A pure spinor $\Phi$ is {\it nondegenerate} if $\ker\Phi\cap\ol{\ker\Phi}=\{0\}$, i.e., 
$(TM\oplus T^*M)^\C=\ker\Phi\oplus\ol{\ker\Phi}$.
Then a nondegenerate, pure spinor $\Phi\in \w^\bullet T^*M$ gives an almost generalized 
complex structure $\J_{\Phi}$ which satisfies 
$$
\J_\Phi E =
\bgn{cases}
&+\sqrt{-1}E, \quad E\in \ker\Phi\\
&-\sqrt{-1}E, \quad E\in \ol{\ker\Phi}
\end{cases}
$$
Conversely, a \complex structure $\J$ arises as $\J_\Phi$ for a nondegenerate, pure spinor $\Phi$ which is unique up to multiplication by
non-zero functions.  
Thus there is a one to one correspondence between almost generalized complex structures and non-degenerate, pure spinors modulo multiplication by non-zero functions.
{\it  The canonical line bundle} $K_{\J_\Phi}$ of $\J_\Phi$ is the complex line bundle generated by the non-degenerate, pure spinor $\Phi$.
$\J_\Phi$ is integrable if and only if $d\Phi=E\cdot\Phi$ for $E\in (TM\oplus T^*M)^\C$. 
The {\it type number} of $\J=\J_\Phi$ is a minimal degree of the differential form $\Phi$, which is allowed to change on a manifold.
\bgn{example}\label{ex:GCS from J} Let $J$ be the ordinary complex structure on a complex manifold $X$. 
Then the complex structure $J^*$ on $T^*X$ is given by 
$(J^*\eta)(v):=\eta(Jv)$ for $\eta\in T^*X, v\in TX$ and we obtain a generalized complex structure  $\J_J$ which is defined by   
$$
\J_J=\bgn{pmatrix}J&0\\0&-J^*
\end{pmatrix},
$$
where the canonical line bundle of $\J_J$ is the ordinary canonical line bundle 
$K_J=\w^{n,0}$ which consists of $n$-forms. Thus we have 
$\text{\rm Type $\J_J =n$}.$
\end{example}
\bgn{example}
Let $\ome$ be a real symplectic structure on a $2n$-manifold $M$. 
Then the interior product $i_v\ome$ of a vector $v$ yields an isomorphism
$\til\ome: TM\to T^*M$ which admits the inverse $\til\ome^{-1}$.
Then a generalized complex structure $\J_\ome$ is defined by 
$$\J_\psi=\bgn{pmatrix}0&-\til\ome^{-1}\\
\til\ome&0
\end{pmatrix},$$ 
Then the canonical line bundle of $\J_\ome$ is generated by 
$$\psi=e^{\sqrt{-1}\ome}=1+\sqrt{-1}\ome +\frac1{2!}(\sqrt{-1}\ome)^2+\cdots+\frac1{n!}(\sqrt{-1}\ome)^n,$$ where the  minimal degree of $\psi$ is $0$. 
Thus  we have
$\text{\rm Type $\J_\psi =0$}$
\end{example}
\bgn{example}[The cation of $b$-fields]\label{ex:b-field}
Let $\J$ be a generalized complex structure which is induced from a 
non-degenerate, pure spinor $\phi$.
A $d$-closed real $2$-form $b$ acts on $\phi$ by 
$e^b\cdot\phi$ which is also a non-degenerate, pure spinor. Thus $e^b\phi$ induces 
a generalized complex structure $\J_b$, which is called the action of $b$-field on $\J$.
The {generalized complex structure} $\J_b$ gives the decomposition $TM\oplus T^*M=L_{\J_b}\oplus \ol{L}_{\J_b}$, where 
$L_{\J_b}=\Ad_{e^b}\circ\J\circ\Ad_{e^{-b}}$ and 
$$\Ad_{e^b}=\bgn{pmatrix}1&0\\b&1.
\end{pmatrix}
$$
\end{example}
\bgn{example}[Poisson deformations]\label{ex:Poisson deformations}
Let $X$ be a complex manifold and
$\b$  a holomorphic Poisson structure on a complex manifold $X$.
Then $\b$  gives deformations of new \complex structures by 
$\J_{\b t}:=\Ad_{e^\b t}\circ \J_J\circ\Ad_{e^{-\b}}$ where 
$$
\Ad_{e^{\b}}=\bgn{pmatrix}1&\b\\0&1\end{pmatrix}
$$
The type number is given by 
Type ${\J_{\b t}}_x=n-2$ rank of $\b_x$ at $x\in M$. 
\end{example}

\section{Deformation theory of generalized complex structures}\label{sec:Deformation theory of GCS}
Let $(M, \J)$ be a generalized complex manifold with  the decomposition 
$(TM\oplus T^*M)^\C=L_\J\oplus \ol L_\J$.
The bundle $\ol L_\J$ is a 
Lie algebroid bundle which yields the 
Lie algebroid complex,
$$
0\to \w^0\ol L_\J\overset{d_L}\to\w^1\ol L_\J\overset{d_L}\to\w^2\ol L_\J\overset{d_L}\to\w^3\ol L_\J\to \cdots
$$
It is known that the Lie algebroid complex is the deformation complex of  generalized complex structures. 
In fact, 
$\e\in \w^2\ol L_\J$ gives a deformed isotropic subbundle 
$L_\e:=\{ E+[\e, E]\, |\, E\in L_\J\}$ which yields a decomposition 
$(TM\oplus T^*M)^\C=L_\e\oplus \ol L_\e$ if $\e$ is sufficiently small.
The isotropic bundle $L_\e$ yields a generalized complex structure if and only if $\e$ satisfies the generalized Mauer-Cartan equation
$$
d_L\e+\frac12[\e, \e]_{\Sch}=0,
$$
where $[\e, \e]_{\Sch}$ denotes the Schouten bracket. 
The Lie algebroid complex $(\w^\bullet, d_L)$ is an elliptic complex and 
the spaces of semi-universal deformations (the Kuranishi spaces) of generalized complex structures are constructed \cite{Gualtieri_2011}.
The space of infinitesimal deformations of {generalized complex structure} is given by the second cohomology group
$H^2(\w^\bullet\ol L_\J)$  and the obstruction spaces of {generalized complex structure} is the third one
$H^3(\w^\bullet\ol L_\J)$.
\bgn{remark}\label{re:Lie algebroid iso b-fields}
Let$\J_b$ be the {generalized complex structure} by the action of $d$-closed $b$ fields on $\J$.
Then the isomorphism $\Ad_{e^b}\,:\,\ol {L}_{\J}\cong \ol {L}_{\J_b}$ yields the isomorphism of the Lie algebroid complexes $\w^\bullet\ol{L}_{\J}\cong \w^\bullet\ol {L}_{\J_b}$. 
Thus we have an isomorphism of the Lie algebroid cohomology groups 
$H^k(\w^\bullet\ol L_\J)\cong H^k(\w^\bullet\ol L_{\J_b})$.
\end{remark}
Let 
$X=(M,J)$ be a complex manifold and $\J:=\J_J$ denotes the generalized complex structure induced from $J$ as in Example \ref{ex:GCS from J}.
Then it turns out that $H^\bullet(\w^\bullet\ol L_\J)$ is the hypercohomology group of the trivial complex of sheaves: 
$$
0\to {\cal O}_X\overset{0}\to \Theta\overset{0}\to\w^2\Theta\overset{0}\to \w^3\Theta\overset{0}\to \cdots,
$$
where $0$ denotes the zero map.
Thus we have 
$$
H^k(\w^\bullet\ol L_{\J_J})=\oplus_{p+q=k}H^p(X, \w^q\Theta),
$$
where $\Theta$ denotes the sheaf of holomorphic vector fields on $X$ and 
$\w^q\Theta$ denotes the $q$-th skew-symmetric tensor of $\Theta$.
The infinitesimal deformations is given by 
$$H^2(\w^\bullet\ol L_{\J_J})=H^2(X, \O_X)\oplus H^1(X, \Theta)\oplus H^0(X, \w^2\Theta),$$
where $H^1(X,\Theta)$ is the infinitesimal deformations of ordinary complex structures and 
$H^2(X, {\cal \O})$ is given by the action of $b$-fields and $H^0(X, \w^2\Theta)$ corresponds to deformations given by holomorphic Poisson structures.\\
A holomorphic $2$-vector 
$\b\in H^0(X, \w^2\Theta)$ is a holomorphic Poisson structure if 
$[\b, \b]_{\Sch}=0$, where $[\, ,\,]_{\Sch}$ stands for the Schouten bracket. A holomorphic $2$-vector 
$\b$ gives the Poisson bracket of functions by
$\{f, g\}_\b=\b(df\w dg)$. Then the Poisson bracket satisfies the Jacobi identity if and only if $\b$ is a Poisson structure.
A holomorphic Poisson structure $\b$ satisfies the generalized Mauer-Cartan equation since $d_L\b=\ol\pa\b=0$. Thus 
$\J_{\b t}=e^{\b t} \J_J e^{-\b t}$ gives deformations of generalized complex structures, where $t$ denotes the complex parameter of deformations.
We denote by $\ol L_{\J_\b}$  the Lie algebroid bundle of $\J_{\b}$.
A holomorphic Poisson structure $\b$ and the Schouten bracket give a map $\del_\b: \w^p\Theta\to \w^{p+1}\Theta$ by 
$\del_\b\a:=[\b,\a]_{\Sch}$, where $\a\in \w^p\Theta$. 
Since the Schouten bracket satisfies the super Jacobi identity and $[\b,\a]_{\Sch}=[\a,\b]_{\Sch}$, we have
$$\del_\b\circ\del_\b(\a)=[\b,[\b,\a]_{\Sch}]_{\Sch}=\frac 12[\a, [\b,\b]_{\Sch}]_{\Sch}=0.
$$
Then we obtain the Poisson complex:
$$
0\to\O_X\overset{\delta_\b}\to\Theta\overset{\delta_\b}\to\w^2\Theta\overset{\delta_\b}\to\w^3\Theta\to\cdots
$$
where $\del_\b f=[\b, f]_{\Sch}=[df,\b]\in \Theta$ for $f\in \O_X$ and
 $[df, \b]$ is the commutator  in the Clifford algebra which is equal to the coupling between $df$ and $\b$. 
 A holomorphic Poisson structure $\b$ defines a map $\til \b$ from the sheaf of holomorphic $1$-forms $\Ome^1$ to 
 $\Theta$ by $[\t,\b]$ for $\t\in \Ome^1$. 
 The map $\til \b$ gives the map $\w^p\til\b : \Ome^p\to \w^p\Theta$ by 
 $$\til\b(\t_1\w\cdots\w\t_p)=\til\b(\t_1)\w\cdots\w\til\b(\t_p).$$

\bgn{proposition}\label{prop:Poisson complex and derham complex}
The map $\w^p\til\b$ induces a map from the holomorphic de Rham complex $(\Ome^\bullet, d)$ to the Poisson complex 
$(\w^\bullet\Theta, \del_\b)$, 
$$
\xymatrix{
0\ar[r]&{\cal O}_X\ar[d]_{\text{\rm id}}\ar[r]^d&\Ome^1\ar[d]_{\til\b}\ar[r]^d&\Ome^2\ar[d]_{\w^2\til\b}\ar[r]^d&\Ome^3\ar[d]_{\w^3\til\b}\ar[r]^d&\cdots\\
0\ar[r]&{\cal O}_X\ar[r]^{\del_{\b}}&\Theta\ar[r]^{\del_{\b}}&
\w^2\Theta\ar[r]^{\del_{\b}}
&\w^3\Theta\ar[r]^{\del_{\b}}
&\cdots
}
$$
\end{proposition}
\bgn{proof}
It follows from our definitions that $\til\b(df)=[df, \b]=[\b, f]_{\Sch}=\del_\b f$ and $\del_\b\circ\del_\b(f)=\del_\b(\til\b(df))=0$ for $f\in {\cal O}_X$.
Since $\del_\b$ is a derivation, we have
\bgn{align}\w^{p+1} \til\b\circ d (f dg_1\w\cdots\w dg_p) =& 
\til\b(df)\w\til\b(dg_1)\w\cdots\w\til\b(dg_p)\\
=&\del_\b( f\til\b(dg_1)\w\cdots\w\til\b(dg_p))\\
=&\del_\b\circ\w^p\til\b(f dg_1\w\cdots\w dg_p),
\end{align}
where $g_1, \cdots, g_p\in {\cal O}_X$. 
Thus we obtain $\w^{p+1}\til\b(d\gam)=\del_\b\circ\w^p\til\b(\gam)$ for $\gam\in \w^p\Theta$.
\end{proof}

The following is already obtained in \cite{Laurent_Stienon_Ping}
\bgn{proposition}\text{\rm\cite{Laurent_Stienon_Ping}}\label{prop:hypercohomology} The hypercohomology of the Poisson complex is isomorphic to the cohomology groups of  the Lie algebroid complex of $\J_\b$, i.e., 
$$\Bbb H^k(\w^\bullet\Theta)\cong H^k(\w^\bullet\ol L_{\J_\b})$$ for all $k$.
\end{proposition}
\bgn{proof}
For the completeness of the paper, we shall give a proof of 
Proposition \ref{prop:hypercohomology}.
Let $\w^{0,1}_\b$ be a vector bundle which is the twisted $\w^{0,1}$ by the adjoint action of $\ol\b$, 
$$
\w^{0,1}_\b:=\{e^{\ol\b} \cdot \ol\t\cdot e^{-\ol\b}=\ol \t+[\ol\b , \ol\t]\in \w^{1,0}\oplus T^{0,1}M\, |\, \ol \t\in \w^{0,1}\,\}.
$$
We denote by $\w^{0,q}_\b$ the $q$-th skew symmetric tensor of $\w^{0,1}_\b$.
By the twisted Dolbeault operator $\ol\pa_\b:=e^{\ol\b}\circ \ol\pa \circ e^{-\ol\b}$, we obtain 
the twisted Dolbeault complex: $(\w^{0,\bullet}_\b, \ol\pa_\b)$. 
Let $T^{1,0}X$ be the tangent bundle on $X$ and $T^{p,0}X$ the skew-symmetric tensor of $T^{1,0}X$. 
Then we also have the complex $(T^{p,0}X\otimes\w^{0,\bullet}_\b, \ol\pa_\b)$ which are a resolution of $\w^p\Theta$. 
As in the Poisson complex, the map $\del_\b$ defines the map between complexes
$\del_\b : (T^{p,0}X\otimes\w^{0,\bullet}_\b, \ol\pa_\b)$ to $(T^{p+1,0}X\otimes\w^{0,\bullet}_\b, \ol\pa_\b)$.
Thus we have a double complex $(T^{\bullet,0}X\otimes \w^{0,\bullet}, \del_\b, \ol\pa_\b)$.
It turns out that the total complex of the double complex is the Lie algebroid complex $(\w^\bullet, \ol L_{\J_\b}, d_L)$. 
Since the cohomology groups of the total complex is the hypercohomology groups of the Poisson complex, we have 
$\Bbb H^k(\w^\bullet\Theta)\cong H^k(\w^\bullet\ol L_{\J_\b})$.

\end{proof}

\section{Log symplectic structure $\ome_{\C}$ in holomorphic category and generalized complex structure $\J_{\phi}$}
Let $X$ be a complex manifold of complex dimension $n$. 
Let $D$ be a smooth divisor on $X$ which admits holomorphic coordinates 
$(z_1, \cdots, z_{2m})$ such that 
$D$ is given by $z_1=0$. We call such coordinates
$(z_1,\cdots, z_{2n})$ logarithmic coordinates.
{\it A logarithmic symplectic structure${}^{\dag}$}\footnote{${}^\dag$The notion of logarithmic symplectic structures was introduced in \cite{Goto_2002}. } $\ome_\C$ is a $d$-closed, logarithmic $2$-form on $X$ which satisfies the followings (1) and (2) :\\
(1) There exist logarithmic coordinates $(z_1, \cdots, z_{2m})$ on a neighborhood of every point in $D$ 
such that
$\ome_\C$ is written as 
\bgn{equation}\label{eq:log.sym}
\ome_\C=\frac{dz_1}{z_1}\w dz_2+dz_3\w dz_4+\cdots+dz_{2m-1}\w dz_{2m}.
\end{equation}
(2) On the complement $X\bsh D$, $\ome_\C$ is a holomorphic symplectic form.\\
Then $\phi:=e^{\ome_\C}$ is a $d$-closed, non-degenerate, pure spinor on $X\bsh D$ which induces 
the generalized complex structure $\J_\phi$ on $X\bsh D$.
On a neighborhood of $D$,  $z_1e^{\ome_\C}$ is also a nondegenerate, pure spinor. 
Hence  $\J_\phi$ can be extended as a generalized complex structure on $X$.
 The type number of $\J_\phi$ is given by the followings: 
On the divisor $D=\{z_1=0\}$, $\J_\phi$ is induced from $z_1\phi$, where
$$
z_1e^{\ome_\C}|_{z_1=0}=dz_1\w dz_2+\cdots=(dz_1\w dz_2)\w e^{\til\ome_\C},
$$
where $\til\ome_\C=dz_3\w dz_4+\cdots+dz_{2m-1}\w dz_{2m}.$
Thus the minimal degree of $z_1\phi$ on $D$ is euqal to $2$. 
On the complement $X\bsh D$, $\J_\phi$ is given by $\phi$ whose the minimal degree is $0$. 
Thus we have
$$\text{\rm Type }\J_{\phi} (x)=\bgn{cases} &2\qquad (x\in D)\\
&0\qquad (x\notin D)
\end{cases}
$$
\section{Lie algebroid cohomology groups of $\J_\phi$ and logarithmic Poisson structure $\b$}
Let $\ome_\C$ be a logarithmic symplectic structure on a complex manifold $X$. 
Then $\ome_\C$ gives the isomorphism between the sheaf of holomorphic logarithmic tangent vectors $\Theta(-\log D)$ and 
the sheaf of logarithmic $1$-forms $\Ome^1(\log D)$ which also indues the isomorphism 
$\w^2\Theta(-\log D)\cong \Ome^2(\log D)$. 
Then $\ome_\C\in H^0(X, \Ome^2(\log D))$ admits the dual $2$-vector 
$\b\in H^0(X, \w^2\Theta(-\log D))$. 
Since $\ome_\C$ is $d$-closed, $\b$ is a Poisson structure. We call the $\b$ {\it a holomorphic log Poisson structure}.

The interior product $i_v\ome_\C$ by a holomorphic vector field $v$ of $\ome_\C$ gives a meromorphic $1$-form
with simple pole along $D$. Thus $\ome_\C$ yields a map $\til\ome_\C$ from $\Theta$ to $\h\Ome^1$,
where $\hat\Ome^1$ is defined to be the image of $\til\ome_\C$. 
Then the map $\til\ome_\C$ is the inverse of the map $\til\b :\h\Ome^1\to \Theta$ as in 
Section \ref{sec:Deformation theory of GCS}.
We define $\hat\Omega^p$ by $\hat\Omega^p=\w^p\hat\Omega^1$. 
Then the exterior derivative $d$ gives a complex $(\hat\Omega^\bullet, d)$. 
Since $\w^p\hat\ome_\C$ gives an isomorphism  $\w^p\Theta\cong\hat\Omega^p$ which is the inverse of $\w^p\til\b$, 
it follows from Proposition \ref{prop:Poisson complex and derham complex} that the complex $(\w^\bullet\hat\Omega, d)$ is isomorphic to the Poisson complex $(\w^\bullet\Theta, \del_\b)$ and we obtain 
$\Bbb H^k(\w^\bullet\Theta)\cong \Bbb H^k(\hat\Ome^\bullet)$.
We denote by $(\Ome^\bullet(\log D), d)$ the holomorphic log complex. 

Since the log complex $(\Ome^\bullet(\log D), d)$ is a subcomplex of $(\h\Ome^\bullet, d)$, we have the short exact sequence of complexes:
$$
0\to \Ome^\bullet(\log D)\to \h\Ome^\bullet\to Q^\bullet\to 0, 
$$
where $Q^\bullet$ denotes the quotient complex. 
\bgn{lemma}\label{lem: Q}
Let ${\cal H}^\bullet(Q^\bullet)$ be the cohomology sheaves of the complex $Q^\bullet$.
Then we have that ${\cal H}^\bullet(Q^\bullet)=\{0\}.$
\end{lemma}
\bgn{proof}
Let $(z_1,\cdots, z_{2n})$ be logarithmic coordinates of a neighborhood of $x\in D$
such that $\ome_\C$ is given by (\ref{eq:log.sym}). 
Then we see that 
$$
i_{\frac{\pa}{\pa z_1}}\ome_\C=\frac{dz_2}{z_1},\qquad  i_{\frac{\pa}{\pa z_2}}\ome_\C=-\frac{dz_1}{z_1}
$$
It follows that the germ of the image $\h\Ome^1_x$ is generated by 
$$\dstyle\lan \frac{dz_2}{z_1}, \frac{dz_1}{z_1}, dz_3,\cdots, dz_{2m}\ran$$ 
over ${\cal O}_{X,x}$.
Then we see that the germ of $\hat\Ome^2_x$ is generated by 
$$\dstyle\lan \frac{dz_1}{z_1}\w\frac{dz_2}{z_1}, \frac{dz_1}{z_1}\w dz_i, \frac{dz_2}{z_1}\w dz_j, dz_i\w dz_j\, |\, i, j=1,\cdots 2m\ran$$ over ${\cal O}_{X,x}$.
Then it turns out that every $\a\in \hat\Ome^2_x$ is written as 
$$
\a=\a_0\w \frac{dz_1}{z_1}\w\frac{dz_2}{z_1}+\a_1\w \frac{dz_1}{z_1}+\a_2\w \frac{dz_2}{z_1}+\a_3
$$
where $\a_0, \a_1,\a_2, \a_3$ are holomorphic forms. 
Since $\dstyle{d(\frac{dz_2}{z_1})=-\frac{dz_1}{z_1}\w\frac{dz_2}{z_1}}$, 
it follows that 
$$
\gam :=\a+(-1)^{|\a_0|}d(\a_0\w\frac{dz_2}{z_1})\in\h\Ome^2_x
$$
is a meromorphic form with pole of order at most one on $D$. 
We assume that $d\a$ is a logarithmic form. Then $z_1d\a$ is holomorphic and 
$z_1d\gam$ is also holomorphic. 
Since $z_1\gam$ is holomorphic, it follow that $\gam$ is a logarithmic form. 
If $\a$ is a representative of the germ of the cohomology sheaves ${\cal H}^\bullet(Q^\bullet)_x$, 
then $[\a]=[\gam]\in {\cal H}^\bullet(Q^\bullet)_x$ and $d\a$ is a logarithmic form. 
Thus it follows that $[\a]=[\gam]=0\in {\cal H}^\bullet(Q^\bullet)_x$ since $\gam $ is a logarithmic form  and $\gam\equiv 0$ in
$Q^\bullet$.

\end{proof}
\bgn{proposition}\label{prop:cohomology haOmega}
The complex $(\w^\bullet\h\Ome, d)$ is quasi-isomorphic to the logarithmic complex 
$(\Ome^\bullet(\log D), d)$. 
Thus the cohomology groups $\Bbb H^k(\w^\bullet \Theta)\cong\Bbb H^k(\w^\bullet\h\Ome)\cong H^k(\w^\bullet\ol L_\phi)$ are given by $H^k(X\bsh D, \C)$.
\end{proposition}
\bgn{proof}
It follows from Lemma \ref{lem: Q} that the map $(\Ome^\bullet(\log D), d)\to (\hat\Ome^\bullet, d)$ is a quasi-isomorphism.
Thus we have $\Bbb H^k(\Ome^\bullet(\log D))\cong\Bbb H^k(\hat\Ome^\bullet)$.
It is known that the hypercohomology groups $\Bbb H^k(\w^\bullet\Ome(\log D))$ of the log complex are 
the cohomology groups of the complement $H^k(X\bsh D, \C)$. 
It follows from Proposition \ref{prop:hypercohomology} that $H^k(\w^\bullet \ol L_\phi)\cong \Bbb H^k(\w^\bullet\Theta)$. 
Thus we obtain $H^k(\w^\bullet\ol L_\phi)\cong H^k(X\bsh D,\C)$.
\end{proof}
\section{Unobstructed deformations of generalized complex structures induced from C$^\infty$ logarithmic symplectic structures}

Let $M$ be a differentiable manifold of real dimension $4m$ and $D$ a submanifold of real codimension $2$.
We assume that there is an open cover $M=\cup_\a U_\a$ such that 
each $U_\a$ is an open set of $\C^{2m}$ with complex coordinates $(z^1_\a,\cdots, z^{2m}_\a)$
and $D$ is locally given by $\{z^1_\a=0\}$ for $U_\a\cap D\neq \emptyset$.
We say $(z^1_\a,\cdots, z^{2m}_\a)$ logarithmic coordinates of $D$.
Note that we do not assume that $M$ is a complex manifold. 
In fact,  defining equations of $D$ satisfies $z_\a^1=e^{f_{\a,\b}}z_\b^1$ on 
$U_\a\cap U_\b$, where $f_{\a,\b}$ is a smooth function on $U_\a\cap U_\b$.

\bgn{definition}\label{def:Cinfty log sym}
{\it A $C^\infty$ logarithmic symplectic structure}${}^\ddag$
\footnote{${}^\ddag$ Note that the notion of $C^\infty$ logarithmic symplectic structures is different from the one of singular symplectic structures as in \cite{{Guillemin_Miranda_Pires}, Gualtieri_Li}
 whose singular loci are real codimension $1$.
}
$\ome_\C$ is a $d$-closed complex $2$-form 
which satisfies the followings:\\
(1) On a neighborhood of $D$, $\ome_\C$ is locally given by  
$$\ome_\C=\frac{dz_1}{z_1}\w dz_2+dz_3\w dz_4+\cdots+dz_{2m-1}\w dz_{2m},$$
where $(z_1,\cdots, z_{2m})$ are logarithmic coordinates of 
$D=\{z_1=0\}$. 
\\
\noindent
(2) On a neighborhood of the complement $M\bsh D$, $\ome_\C=b+\sqrt{-1}\ome$  
where $b$ is a $d$-closed $2$-form and $\ome$ denotes a real symplectic structure.
\end{definition}
Then $\phi =e^{\ome_{\C}}$ is a $d$-closed, non-degenerate, pure spinor which induces the generalized complex structure $\J_\phi$ on the complement $M\bsh D$.
In fact, $z_1\phi$ is a non-degenerate pure spinor on a neighborhood $U$ of $D$. 
Thus it follows that $\phi:=e^{\ome_\C}$ defines a generalized complex structure $\J_\phi$ on $M$.

Then we have the following theorem :
\bgn{theorem}\label{th:log.cohomology}
Let $\ome_\C$ be a $C^\infty$ logarithmic symplectic structure on $M$ and $\J_\phi$ the generalized complex structure which is induced from
$\phi=e^{\ome_\C}$. 
Then the Lie algebroid cohomology groups $H^k(\w^\bullet\ol L_\phi)$ of $\J_\phi$ is isomorphic to 
$H^k(M\bsh D, \C)$.
\end{theorem}
\bgn{theorem}\label{th:unobst. log GCS}
Let $\ome_\C$ be a $C^\infty$ logarithmic symplectic structure on $M$ and $\J_\phi$ the  generalized complex structure which is induced from
$\phi=e^{\ome_\C}$. 
Then
deformations of $\J_{\phi}$ are unobstructed and the space of infinitesimal deformations is given by 
$H^2(M\bsh D, \C)$.
\end{theorem}

In order to prove our theorems, we shall introduce  $C^\infty$ logarithmic deformations of $\J_\phi$ and
the $C^\infty$ logarithmic deformations are unobstructed in next Section.
\section{Logarithmic deformations of $\J_\phi$}
A $C^\infty$ logarithmic vector field $V$ on a manifold $M$ along $D$ is a $C^\infty$ vector field which is locally given by 
$$
V=f_1z_1\frac{\pa}{\pa z_1}+g_1\frac{\pa}{\pa \ol{z_1}}+\sum_{i=2}^{2m}
f_i\frac{\pa}{\pa z_i}+g_i\frac{\pa}{\pa \ol z_{i}},
$$
where $f_i,g_i\, (i=1,\cdots, 2m)$ are $C^\infty$ functions and $(z_1, \cdots, z_{2m})$ are logarithmic coordinates 
of $D=\{z_1=0\}$.
Thus a $C^\infty$ logarithmic vector field $V$ preserves the ideal $(z_1)$ which is an analog of the notion of
logarithmic vector fields in complex geometry.
We denote by ${\cal T}_{\log}M$ the sheaf of $C^\infty$ logarithmic vector fields. 
The sheaf ${\cal T}_{\log}M$ is locally free which gives a $C^\infty$ vector bundle $T_{\log}M$.
Our generalized complex structure $\J_\phi$ gives the decomposition 
$TM\oplus T^*M=L_\phi\oplus \ol L_\phi$ and the Lie algebroid complex 
$(\w^\bullet \ol L_\phi, d_{L_\phi})$.
We define a subbundle $\ol L^{\log}_\phi$ by 
$$
\ol L_\phi^{\log}=\ol L_\phi\cap (T_{\log}M\oplus T^*M)^\C
$$
Then we obtain the subcomplex $(\w^\bullet\ol L_\phi^{\log}, d_{L})$ of the Lie algebroid complex $(\w^\bullet \ol L_\phi, d_L)$.
Then we have 
\bgn{proposition}\label{prop:log lielagebroid cohomology}
The cohomology group $H^k(\w^\bullet\ol L_\phi^{\log})$ of the subcomplex $(\w^\bullet\ol L_\phi^{\log}, d_{L})$ is isomorphic to the cohomology group $H^k(M\bsh D, \C)$
\end{proposition}
\bgn{proof}Let ${\cal \ol L}_\phi^{\log}$ be the sheaf of $C^\infty$ sections of the bundle $\ol L_\phi^{\log}$
and $\w^p{\cal \ol L}_\phi^{\log}$ the $p$-th skew symmetric tensor of ${\cal \ol L}_\phi^{\log}$. 
Then we have the complex of sheaves: 
\bgn{equation}\label{eq:sheaves log Lie algebroid complex}
0\arrow \w^0{\cal \ol L}_\phi^{\log}\overset{d_L}\arrow \w^1{\cal \ol L}_\phi^{\log}\overset{d_L}\arrow \w^2{\cal \ol L}_\phi^{\log}
\overset{d_L}\arrow\cdots. 
\end{equation}
Since $\w^p{\cal \ol L}_\phi^{\log}$ is a soft sheaf, the hypercohomology groups 
$\Bbb H^k(\w^\bullet{\cal \ol L}_\phi^{\log})$ of the complex of sheaves $(\w^\bullet{\cal \ol L}_\phi^{\log}, d_L)$
is given by global sections and we have $\Bbb H^k(\w^\bullet{\cal \ol L}_\phi^{\log})\cong H^k(\w^\bullet\ol L_\phi^{\log})$.
The interior product $i_v\ome_{\C}$ of a vector $v$ by $\ome_\C$ restricted  to the complement $M\bsh D$ gives a map from $T_{\log}M$ to $1$-forms on $M\bsh D$ which induces a map $\til{\ome_{\C}}$ from ${\cal \ol L}_\phi^{\log}$ to the sheaf of differential $1$-forms on $M\bsh D$ by 
$\til{\ome_\C}(v+\t)=-i_v{\ome_\C}+\t.$ Then we have the map 
$\w^p\til\ome_\C : \w^p{\cal \ol L}_\phi^{\log}\to \A^p(M\bsh D)$ by 
\bgn{equation}\label{eq:wp til omeC}
\w^p\til\ome_\C(v_1\w\cdots\w v_s\w\t_1\w\cdots\w\t_s) =(-1)^si_{v_1}\ome_\C\w\cdots\w i_{v_s}\ome_\C\w\t_1\w\cdots\w\t_s,
\end{equation}
where $s+t=p$ and $\A^p(M\bsh D)$ denotes the sheaf of $p$-forms on $M\bsh D$.
The map $\w^p\til\ome_\C$ gives the map $\w^\bullet\til\ome_\C$ from the complex
$(\w^\bullet{\cal \ol L}_\phi^{\log}, d_L)$
to the de Rham complex $(\A^\bullet(M\bsh D), d)$.
We shall show that the map 
$\w^\bullet\til\ome_\C :(\w^\bullet{\cal \ol L}_\phi^{\log}, d_L)\to (\A^\bullet(M\bsh D), d)$ is quasi-isomorphic. 
In order to obtain the quasi-isomorphism, we shall determine the cohomology sheaves ${\cal H}^k(\w^\bullet {\cal L}_\phi^{\log})(U)$ of the complex $(\w^\bullet{\cal \ol L}_\phi^{\log}, d_L)$ restricted to a neighborhood $U$ in the followings two cases :
\\
(1) If $U$ is a neighborhood of $D$ with logarithmic coordinates $(z_1,,\cdots, z_{2m})$, then 
the logarithmic coordinates define the complex structure on $U$ such that $\ome_\C|U$  is a
logarithmic symplectic structure which is the dual of holomorphic logarithmic Poisson structure $\b$ as in Section 5. 
Then it turns out that the cohomology groups ${\cal H}^k(\w^\bullet {\cal L}_\phi^{\log})(U)$ is given by 
the hypercohomology groups $\Bbb H^k(\w^\bullet(\Theta(-\log D))$ of the Poisson complex of multi-logarithmic tangent vectors: 
 $$
 0\to \w^0\Theta(-\log D)\overset{\del_\b}\arrow\w^1\Theta(-\log D)\overset{\del_\b}\arrow\w^2\Theta(-\log D)\overset{\del_\b}\arrow\cdots
 $$
 The map $\w^\bullet\til\ome_\C$ restricted to $U$ gives an isomorphism from the complex of logarithmic multi-tangent vectors
 $(\w^\bullet\Theta(-\log D),\del_\b)$ to the complex of logarithmic forms $(\Ome^\bullet(\log D),d)$ which 
 induces the isomorphism between cohomology groups $\Bbb H^k(\Ome^\bullet(\log D))\cong \Bbb H^k(\w^\bullet\Theta(-\log D)).$
 It is known that the hypercohomology groups $\Bbb H^k(\Ome^\bullet(\log D))$ of the complex of logarithmic forms 
is $H^k(U\bsh U\cap D, \C)$. 
Thus we have ${\cal H}^k(\w^\bullet {\cal L}_\phi^{\log})(U)\cong H^k(U\bsh U\cap D, \C)$.
\\
(2) If $U$ is a neighborhood of the complement $M\bsh D$, then $\ome_\C|_U$ is given by $b+\sqrt{-1}\ome$.
It follows that the map $\w^\bullet\til\ome_\C: (\w^\bullet{\cal \ol L}_\phi^{\log}(U), d_L)\to 
(\A^\bullet(U), d)$ is an isomorphism and ${\cal H}^k(\w^\bullet {\cal L}_\phi^{\log})(U)\cong H^k(U, \C)$.

Then it follows that the map $\w^\bullet\til\ome_\C :(\w^\bullet{\cal \ol L}_\phi^{\log}, d_L)\to (\A^\bullet(M\bsh D), d)$ is a quasi-isomorphism. 
Hence we obtain $H^k(\w^\bullet\ol L_\phi^{\log})\cong H^k(M\bsh D, \C)$.
\end{proof}

\bgn{proposition}\label{prop:unobstructed def of log}
The second cohomology group $H^2(\w^\bullet\ol L_\phi^{\log})\cong H^2(M\bsh D,\C)$ gives unobstructed deformations of generalized complex structures.
\end{proposition}
\bgn{proof}
A $C^\infty$ logarithmic $1$-form $\t$ is a $C^\infty$ $1$-form on $M\bsh D$ which is written on a neighborhood $U$ of $D$ by 
$$
\t=f_1\frac{dz_1}{z_1}+g_1 d\ol z_1+\sum_{i=2}^{2m} f_i dz_i+ g_i d\ol z_i,
$$
where $f_i, g_i$ are $C^\infty$ functions and $(z_1, \cdots, z_{2m})$ denotes logarithmic coordinates. 
Let ${\cal T}^*_{\log}M$ be the sheaf of $C^\infty$ sections of $C^\infty$ logarithmic $1$-forms and $\w^pT^*_{\log}M$ the $p$-th skew symmetric tensors of ${\cal T}^*_{\log}M$. 
Then we have the complex of $C^\infty$ logarithmic forms:
$$
0\to \w^0{\cal T}^*_{\log}M\overset{d}\arrow \w^1{\cal T}^*_{\log}M\overset{d}\arrow \w^2{\cal T}^*_{\log}M\overset{d}\arrow \cdots
$$
It turns out that the hypercohomology groups $\Bbb H^k(\w^\bullet {\cal T}^*_{\log}M)$ of the complex of $C^\infty$ logarithmic forms
is $H^k(M\bsh D, \C)$.
Thus every element of $H^2(M\bsh D, \C)$ admits a $d$-closed representative 
$\a\in C^\infty(M, \w^2{\cal T}^*_{\log}M)$. 
If $\a$ is sufficiently small, then $e^{\ome_\C+ \a }$ is a $d$-closed, non-degenerate, pure spinor on $M\bsh D$ which induces a family of deformations of generalized complex structures  on $M\bsh D$ parametrized by an open set of $H^2(M\bsh D, \C)$.
On a neighborhood of $D$ with logarithmic coordinates $(z_1, \cdots, z_{2m})$, 
$\ome_\C$ is given by 
$\dstyle{\ome_\C=\frac{dz_1}{z_1}\w dz_2+\hat\ome_\C}$ and
$\a$ is written as 
$\dstyle{\a=\frac{dz_1}{z_1}\w\a_1+\gam}$, where $\hat\ome_\C, \a_1,\gam$ are $C^\infty$ $2$-forms. 
Then we have $\dstyle{\ome_\C+\a=\frac{dz_1}{z_1}\w(dz_2+\a_1)+\hat\ome_\C+\gam}.$
Then $z_1e^{\ome_\C+\a}$ restricted to $D=\{z_1=0\}$ is given by 
$$
z_1e^{\ome_\C+\a}|_{z_1=0}=dz_1\w(dz_2+\a_1)\w e^{\hat{\ome_\C}+\gam}
$$
Hence $z_1e^{\ome_\C+\a}|_{z_1=0}$ is also a non-degenerate, pure spinor on $D$ for 
sufficiently small $\a$.
Thus $e^{\ome_\C+ \a }$ gives deformations of generalized complex structures on $M$ which are parametrized by an open set $H^2(M\bsh D,\C)$.
\end{proof}
\section{Proof of main theorems}

\bgn{proof}[\indent\sc Proof of Theorem \ref{th:log.cohomology}]
Let $\ol{\cal L}_\phi$ be  the sheaf of germs of smooth sections of the bundle $\ol L_\phi$. 
Then the Lie algebroid complex gives the complex of sheaves:
\bgn{equation}\label{sheaves of Liealgebroid complex}
0\to \w^0{\cal \ol L_\phi}\overset{d_L}\to\w^1{\cal \ol L_\phi}\overset{d_L}
\to\w^2{\cal \ol L_\phi}\overset{d_L}\to\cdots 
\end{equation}
The hypercohomology groups of the complex of sheaves are isomorphic to the cohomology groups of the Lie algebroid complex since $\w^p{\cal \ol L}_\phi$ are a soft sheaf.
We shall apply the similar argument as in Proposition \ref{prop:log lielagebroid cohomology} to the complex (\ref{sheaves of Liealgebroid complex}).
The interior product $i_v\ome_\C$ of a vector field $v$ by $\ome_\C$ restricted to $M\bsh D$ gives the map $\w^p\til\ome_\C :\w^p{\cal \ol L}_\phi \to 
\A^p(M\bsh D)$ as in (\ref{eq:wp til omeC}) which yields a map $\w^\bullet\til\ome_\C$
form the complex of sheaves $(\w^\bullet{\cal \ol L}_\phi, d_L)$ to 
 the de Rham complex $(\A^\bullet(M\bsh D), d)$.
 We shall  show the map $\w^\bullet\til\ome_\C$ is a quasi-isomorphism. 
 
The sheaves ${\cal H}^\bullet(\w^\bullet{\cal \ol L}_\J)$ of cohomology groups of (\ref{sheaves of Liealgebroid complex}) are determined by the following two cases (1) and (2): 
\\
(1) If $U$ is a neighborhood of $D$ admitting logarithmic coordinates of $D$, then 
the logarithmic coordinates define the complex structure on $U$ such that $\ome_\C|U$  is a
logarithmic symplectic structure which is the dual of holomorphic logarithmic Poisson structure $\b$ as in Section 5. 
It follows from Proposition \ref{prop:cohomology haOmega} that the cohomology 
${\cal H}^k(\w^\bullet{\cal \ol L}_\J)(U)$ is isomorphic to $H^k(U\bsh U\cap D, \C)$.
Thus
if $x\in M\bsh D$, then $\w^\bullet\til\ome_\C$ induces an isomorphism
$${\cal H}^k(\w^\bullet{\cal L}_\J)_x\cong \lim_{\overset{\arrow} {x\in U}} H^k(U) 
$$
\\
\noindent
(2) If $U$ is a neighborhood of the complement $M\bsh D$, then $\ome_\C=b+\sqrt{-1}\ome$ gives an isomorphism from the complex 
$(\w^\bullet{\cal \ol L}_\phi(U), d_L)$ to the de rham complex $(\A^\bullet(U), d)$.
\\
It follows that the complex of sheaves $\w^\bullet{\cal L}_\J$ is quasi-isomorphic to 
the complex of sheaves $(\A^\bullet(M\bsh D), d)$. 
Thus the hypercohomology groups of the complex $\w^\bullet{\cal L}_\J$ are 
$H^\bullet(M\bsh D)$. 
Hence the cohomology groups of the Lie algebroid complex are also $H^\bullet(M\bsh D)$.  
\end{proof}
\bgn{proof}[\indent\sc Proof of Theorem \ref{th:unobst. log GCS}]
From Theorem \ref{th:log.cohomology}, we obtain $H^2(\w^\bullet\ol L_\phi)\cong H^2(M\bsh D, \C)$. 
It follows from Proposition \ref{prop:unobstructed def of log} that we already obtain deformations of generalized complex structures on $M$ parametrized by an open set $H^2(M\bsh D, \C)$. 
Thus we obtain unobstructed deformations of generalized complex structures.
\end{proof}

\section{Generalized complex structures on $4$-manifolds}
\subsection{Non-degenerate, pure spinors of even type on $4$-manifolds}
Let $M$ be a $4$-dimensional manifold.
Then the spinor inner metric of even forms is defined by 
$$
\lan \phi, \psi\ran:=\phi_0\psi_4-\phi_2\w\psi_2+\phi_4\psi_0\in \w^4T^*M,
$$
where $\phi=\phi_0+\phi_2+\phi_4, \psi=\psi_0+\psi_2+\psi_4\in \w^*T^*M$ and 
$\phi_i,\psi_i\in \w^iT^*M$.
Then the spinor inner metric gives a simple description of
 non-degenerate, pure spinors of even type, that is, 
 $\phi\in \w^{\even}T^*M$ is a non-degenerate, pure spinor if and only if $\phi$ satisfies the followings:
\bgn{align}\label{ali:non-deg,pure spinor on 4-mfd1}
&\lan\phi, \phi\ran:=2\phi_0\w\phi_4-\phi_2\w\phi_2=0, \\
&\lan\phi, \ol{\phi}\ran:=\phi_0\w\ol{\phi_4}+\phi_4\w\ol{\phi_0}-\phi_2\w\ol{\phi_2}\neq0
\label{ali:non-deg,pure spinor on 4-mfd2}
\end{align}
A non-degenerate, pure spinor $\phi$ gives a {generalized complex structure} if and only if 
\bgn{equation}\label{eq:GCS on 4-mfd}
d\phi=e\cdot\phi,
\end{equation} where 
$e=v+\t\in TM\oplus T^*M$.
\bgn{example}
Let $(z_1,z_2)$ be the coordinates of $\C^2$ and $\dstyle{\phi=1+\frac{dz_1}{z_1}\w dz_2}$.  Then
$z_1\phi=z_1+dz_1\w dz_2 $ is a  non-degenerate, pure spinor which induces the generalized complex structure  $\,\J_\phi$.
On the divisor $D:=\{z_1=0\}$, we have $\phi|_{D}=dz_1\w dz_1$ and $\J_\phi|_D$ is  induced from a complex structure and  Type $\J_\phi|D=2.$ 
On the complement $\C^2\bsh D$, 
$\dstyle{\phi=\exp(\frac{dz_1}{z_1}\w dz_2)=\exp(b+\sqrt{-1}\ome)}$ , where 
$b$ is a $d$-closed real $2$-form and $\ome$ is a symplectic form. Thus  $\J_\phi$ is coming from a symplectic structure twisted by the action of  $b$-field and
Type $J_\phi|_{\C^2\bsh D}=0$
\end{example}

\subsection{Unobstructed deformations of generalized complex structures on Poisson surfaces}
Let $S=(M, J)$ be a complex surface with effective anti-canonical line bundle $K^{-1}$ and 
$\J_J$ the {generalized complex structure} given by the complex structure $J$ as in Example \ref{ex:GCS from J}.
Then a section $\b\in K^{-1}$ is a holomorphic Poisson structure on $S$ which gives Poisson deformations $\{\J_{\b t}\}$ of {generalized complex structure}s as in Example \ref{ex:Poisson deformations}. 
Then we have 
\bgn{theorem}\label{th:ubobstructed defs of Jb}
If the zero set of $\b$ is a smooth divisor $D$ of $S=(M, J)$, then the Lie algebroid cohomology groups $H^k(\w^\bullet\ol L_{\J_\b})$ is isomorphic to $H^k(M\bsh D, \C)$ and 
deformations of {generalized complex structure} $\J_\b$ are parametrized by an open set $H^2(M\bsh D, \C)$.
\end{theorem}
\bgn{proof}
Let $\ome_\C$ be the dual of $\b$ which is a logarithmic symplectic structure on $S$ of $D$.
Since $\J_\b$ is induced from $\phi=e^{\ome_\C}$,
the results follows from Theorem \ref{th:log.cohomology} and Theorem \ref{th:unobst. log GCS}.
\end{proof}
In general, if an anti-canonical divisor $D$ on a complex surface $S$ is not smooth, the Lie algebroid cohomology groups
(the hypercohomology groups of Poisson complex) are different from the singular cohomology groups of the complement $M\bsh D$.
We assume that $D$ is given by the zero locus of a Poisson structure $\b$ which has isolated singular points $\{p_i\}_{i=1}^m$ and
$\b$ is written as $f_i\frac{\pa}{\pa z_1}\w\frac{\pa}{\pa z_2}$ at a neighborhood of each $p_i$, where $f_i\in {{\cal O}_{S. p_i}}$.
Then we have the complex $(\h\Ome^\bullet, d)$ which is isomorphic to the Poisson complex $(\w^\bullet\Theta, \del_\b)$ as in Section 5.
The log complex $\Ome^\bullet(\log D)$ is a subcomplex of $\h\Ome^\bullet$ and we have the short exact sequence of complexes
\bgn{equation}\label{eq:hOmega to Q}
0\to\Ome^\bullet(\log D)\to \h\Ome^\bullet \to Q^\bullet\to 0
\end{equation}
The cohomology sheaves of the complex $Q^\bullet$ are given by
$$
{\cal H}^i(Q^\bullet)=\bgn{cases}& J_{p_i}\quad (i=2)\\&0\quad \,\,\,\,\,(i\neq 2),  \end{cases}
$$
where $J_{p_i}$ denotes the quotient ring ${\cal O}_{S, p_i}/
(f_{i}, \frac{\pa f_{i}}{\pa z_1}, \frac{\pa f_{i}}{\pa z_2})$ and $(f_{i}, \frac{\pa f_{i}}{\pa z_1}, \frac{\pa f_{i}}{\pa z_2})$ is the ideal generated by $f_i$ and the partial derivatives of $f_i$.
Thus we have 
\bgn{equation}{\Bbb H}^i(Q^\bullet)=\bgn{cases}&\oplus_i J_{p_i}\quad (i=2)\\&0 \qquad\quad (i\neq 2)\end{cases}
\end{equation}
In particular, $p_i$ is a node, then $\J_{p_i}=\C$.
If an anti-canonical divisor $D$ is a simple normal crossing divisor (nodes), then it is known that 
the hypercohomology groups ${\Bbb H}^i(\Ome^\bullet(\log D))$ of the logarithmic complex is given by the singular cohomology groups $H^i(S\bsh D)$ of the complement of the divisor. 
Thus we can calculate the Poisson cohomology groups by using the cohomology groups of the complement and each quotient ring $J_{p_i}$. 
In fact, the short exact sequence (\ref{eq:hOmega to Q}) yields the long exact sequence of the hypercohomology groups.
Thus we have 
$$
{\Bbb H}^i(\h\Ome^\bullet)\cong H^i(S\bsh D)\quad (i=0,1,4).
$$

If $\Bbb H^3(\Ome^\bullet(\log D))\cong H^3(S\bsh D)=\{0\}$, then we have 
$$
0\to H^2(S\bsh D)\to H^2(\h\Ome^\bullet)\to \oplus_i J_{p_i}\to 0
$$
\bgn{proposition}
Let $S$ be a complex surface with a Poisson structure $\b$.
We assume that the divisor $D=\{\b=0\}$ has $m$ nodes $\{p_i\}_{i=1}^m$ and 
$H^3(S\bsh D)=\{0\}$. Then three cohomology groups $\Bbb H^i(\h\Ome^\bullet)\cong\Bbb H^i(\w^\bullet\Theta, \del_\b)\cong H^i(\w^\bullet\ol L_{\J_\b})$ are given by  
$$
\dim H^i(\w^\bullet\ol L_\b)=\bgn{cases}&\dim H^2(S\bsh D)+m\quad (i=2)\\ &\dim H^i(S\bsh D)\qquad\quad (i\neq 2)\end{cases}
$$

\end{proposition}
\bgn{remark}\label{re:cubic curve with one node}
If $D$ is a cubic curve in $\C P^2$ with only a node, then $H^2(\w^\bullet\ol L_{\J_\b})=\C^2$, however $H^2(\C P^2\bsh D)=\C$.
Thus the Lie algebroid cohomology groups of the singular $D$ are different from the cohomology groups of the complement $M\bsh D$.
\end{remark}

\bgn{example}[Del Pezzo surfaces]\label{ex:del Pezzo}
Let $S_k=(M,J)$ be a del pezzo surface which is a blown up $\C P^2$ at generic $k$ points with ample anti-canonical line bundle, where $0\leq k \leq 8$.
If a Poisson structure $\b$ has a smooth divisor $D$, then we have $\dim H^2(S_k\bsh D)=2+k$. 
Thus we have 
\bgn{proposition}
If $D$ is smooth, then $\J_\b$ has $2+k$ dimensional unobstructed deformations.
\end{proposition}
The Poisson cohomology of $S_k$ was already calculated in \cite{Hong_Xu}. 
Their calculation is different from ours. We shall follow their method.
 Let $\b$ be a holomorphic Poisson structure on $S$ which is a section of $K^{-1}$. 
Since the Lie algebroid cohomology groups coincide with the hypercohomology groups of the Poisson complex : 
$$
0\to {\cal O_S}\overset{\del_\b}\to \Theta_S\overset{\del_\b}\to \w^2\Theta\to 0
$$
The hypercohomology groups is the total cohomology groups of the double complex 
$(\ch{C}^p(\w^q\Theta), \del_b, \del)$.  Then the $E_1$-terms are given by the $\ch {C}$ech cohomology groups :
$$
\bgn{matrix} 
H^{2}(S,{\cal O}_S)&H^{2}(S, \Theta)&H^{2}(S, \w^2\Theta)\\
H^{1}(S,{\cal O}_S)& H^{1}(S, \Theta)&H^{1}(S, \w^2\Theta)\\
H^{0}(S,{\cal O}_S)&H^{0}(S, \Theta)&H^{0}(S, \w^2\Theta)
\end{matrix}
$$
Applying the Kodaira vanishing theorem, we obtain that $E_1$-terms are given by
$$
\bgn{matrix}
0&0&0\\
0&H^{1}(S, \Theta)&0\\
\C&H^{0}(S, \Theta)& H^0(S, \w^2\Theta)
\end{matrix}
$$
The map $\del_\b : H^0(S, \Theta)\to H^0(S, \w^2\Theta)$ yields the $E_2$-terms, that is, 
$$
\bgn{matrix}
0&0&0\\
0&H^1(S, \Theta)&0\\
\C&\ker\del_\b&\text{\rm image }\del_\b
\end{matrix}
$$
Thus it follows that the double complex degenerates at $E_2$-terms. 
A holomorphic vector field $\a\in\ker\del_\b$ is given by $\del_\b\a=-{\cal L}_{\a}\b=0$. 
Thus $\a\in \ker\del_\b$ gives an automorphism of $S_k$ preserving a section $\b\in K^{-1}$.
The followings are known for $S_k$ 
$$
\dim H^0(S_k, \Theta)=\bgn{cases}&8-2k\quad (k=0,1,2,3)\\& 0\,\qquad\quad (k=4,5,6,7,8)\end{cases}
$$
$$
\dim H^1(S_k, \Theta)=\bgn{cases}&2k-8\quad (k=5,6,7,8)\\& 0\,\,\quad\qquad (k\leq 4)\end{cases}
$$
$$
\dim H^0(S_k, K^{-1}) =10-k
$$
Thus we obtain $\Bbb H^i(\w^\bullet\Theta)\cong H^i(\w^\bullet\ol L_{\J_\b})$ by 
$$
\dim H^i(\w^\bullet\ol L_{\J_\b})=
\bgn{cases}
&1\quad(i=0)\\
&{\dim\ker\del_\b}\quad(i=1)\\
&\dim H^1(S, \Theta)+\text{\rm image }\del_\b\quad (k=2)\\
&0\quad (k\neq0,1,2)
\end{cases}
$$
If no automorphism preserves an anti-canonical divisor $D$, then $\dim\ker\del_\b=0$. 
In particular, if $k=4,5,6,7,8$, then $\dim H^i(\w^\bullet\ol L_{\J_\b})$ does not depend on a choice of $\b$.
\end{example}

\bgn{remark}
The calculations as in Example \ref{ex:del Pezzo} hold for degenerate del pezzo surfaces which are blown up $\C P^2$ at the set of points in
almost general position. Then $K^{-1}$ is not ample, however we have $H^1(K^{-1})=\{0\}$ (c.f. \cite{Goto_2012}).
\end{remark}
\bgn{example}[Hirzebruch surfaces]
Let $F_e$ be the projective space bundle $\Bbb P({\cal O}\oplus {\cal O}(-e))$ over $\C P^1$ with $e>0$, which is called Hirzebruch surface. Let $f$ be a fibre of $F_e$. Then $K^{-1}$ is given by 
$2b+(e+2)f$, where $b$ is a section of $F_e$ with $b\cdot b=-e$.
Since $K^{-1}$ is effective, we have the $E_1$-terms: 
$$
\bgn{matrix}
0&0&0\\
0&H^1(F_e, \Theta)&H^1(F_e, K^{-1})\\
\C&H^0(F_e, \Theta)&H^0(F_e, K^{-1})
\end{matrix}
$$
These are given by 
\bgn{align*}
&\dim H^1(F_e, \Theta)=e-1, \,\,\dim H^1(F_e, K^{-1})=e-3\\
&\dim H^0(F_e, \Theta)=e+5, \,\,\dim H^0(F_e, K^{-1})=e+6,
\end{align*}
where $\dim H^1(F_e, K^{-1})=0$ if $e\leq 3$. 
Let $\b$ be a poisson structure of $F_e$ which gives a smooth anti-canonical divisor $D$. 
Then it turns out that 
$$
H^i(F_e\bsh D) =\bgn{cases}&\C\quad(i=0)\\&\C^3\quad (i=2)\\
&0\quad (i\neq 0,2) \end{cases}
$$
\bgn{proposition}
Thus we obtain unobstructed deformations of $\J_\b$ which are parametrized by an open set of $H^2(F_e\bsh D)\cong \C^3.$
\end{proposition}
Since the double complex degenerates at $E_2$-terms, it follows from Theorem \ref{th:log.cohomology}
and $H^1(F_e\bsh D)=H^3(F_e\bsh D)=\{0\}$ that the map between $E_1$-terms
$\del_\b^0 : H^0(F_e, \Theta)\to H^0(F_e, K^{-1})$ is injective and the map $\del^1_\b: H^1(F_e, \Theta)\to H^1(F_e, K^{-1})$ is surjective. 
\bgn{remark}
Compared with unobstructed deformations of $\J_\b$, 
it is remarkable  that deformations of {generalized complex structure} starting from $\J_J$ are always obstructed if $e>3$. 
In fact, infinitesimal deformations of $\J_J$ are  $H^1(F_e, \Theta)\oplus H^0(F_e, K^{-1})$ and 
the Kuranish map is given by the Schouten bracket, that is, 
$\a+\b \mapsto [\a+\b, \a+\b]_{\Sch}$, where $\a\in H^1(F_e, \Theta)$ and $\b\in H^0(F_e, K^{-1})$.
Since $H^2(F_e, \Theta)=\{0\}$, the obstruction is $2[\a,\b\]_{\Sch}\in H^1(F_e, K^{-1})$ which is equal to $2\del_\b^1(\a)$. Since the map $\del_\b^1: H^1(F_e, \Theta)\to H^1(F_e, K^{-1})$ is surjective, there exists a $\a\in H^1(F_e, \Theta)$ such that $\del_\b(\a)\neq0$. 
Thus we have 
\bgn{proposition}
There always exists an obstruction to deformations of the generalized complex structure $\J_J$ induced by the ordinary complex structures $J$ on a Hirzebruch surface $F_e$ if $e>3$.
\end{proposition}

\end{remark}
\end{example}

\subsection{$C^\infty$ logarithmic transformations and unobstructed deformations of generalized complex structures}
The natural projection $\w^\bullet T^*M\to \w^0T^*M$ gives rise to a section $s\in \Gam(K^*_\J)$ on $(M, \J)$. Then the zero set
$s^{-1}(0)$ is called  {\it Type changing loci of} $(M, \J)$
which is the set of points of $M$ where the Type number of $\J$ changes from $0$ to $2$ if $\J$ is a {generalized complex structure} of even type on a $4$-manifold.
We denote by $\k(\J)$ {\it the number of connected components of Type changing loci of} $(M, \J).$
Let $(M, \J)$  be a generalized complex $4$-manifold. 
A sub $2$-torus $T\subset M$ is {\it a symplectic torus} if 
there is a neighborhood $\nu(T)$ such that the $\J|{\nu(T)}$ is induced from a real symplectic structure $\ome$ and a $b$-field and 
$T\subset (\nu(T), \ome)$ is a symplectic submanifold.
Let $T\subset (M, \J)$ be a symplectic torus with trivial normal bundle. 
Then
a neighborhood $\nu(T)$ is diffeomorphic to $ D^2\times T^2$ with  
the boundary $\pa\nu(T)\cong T^3$. 
A $C^\infty$ {\it logarithmic transformation} on $T$ consists of two operations:  removing $\nu(T)$ and attaching the manifold $D^2\times T^2$
by a diffeomorphism $\Psi: \pa D^2\times T^2\to \pa\nu(T)$. 
A $C^\infty$ logarithmic transformation on $T$ yields a manifold
$$
M_\Psi:=(M\bsh\,\text{\rm Int}\nu(T))\cup_\Psi (D^2\times T^2)
$$
where Int$\nu(T)$ denotes the interior of $\nu(T)$.
If $\Psi(\pa D^2\times{pt})\subset D^2\times T^2$ is null-homotopic in $\nu(T)$, the $C^\infty$ logarithmic transformation is {\it trivial}. 
The boundary $\pa\nu(T)$ is identified with $\pa D^2\times T^2$ and we denote by $\pi$ the projection
from $\pa D^2\times T^2$ to $\pa D^2$. Then we have the diagram 
$$
\xymatrix{ \pa D^2\times T^2\ar[r]^\Psi\ar[rd] &\pa \nu(T^3)\ar[d]_\pi\\
 &\pa D^2 &}
$$
and the map 
$$
\pi\circ \Psi |_{\pa D^2\times \{pt\}}: \pa D^2\times \{pt\}\to \pa D^2.
$$
Then we define 
 {\it multiplicity} of $C^\infty$ logarithmic transformation to be the  degree of the map $\pi\circ \Psi |_{\pa D^2\times \{pt\}}: S^1\to S^1$

\bgn{theorem}{\rm \cite{Goto_Hayano}}\label{th:Goto and Hayano 1}
Let $(M, \J)$ be a generalized complex $4$-manifold and $T$ a symplectic torus of $(M,\J)$ with trivial normal bundle.
Then every nontrivial $C^\infty$ logarithmic transformation yields  
 a twisted generalized complex structure $\J_\Psi$ on the manifold $M_\Psi$  with $\k(\J_\Psi)=\k(\J)+1$.
 In particular, if $H^3(M_\Psi)=\{0\}$, then $\J_\Psi$ is a generalized complex structure.
\end{theorem}
The proof of this theorem is already given in {\rm \cite{Goto_Hayano}}.
For the completeness of the paper, we will give the proof of the theorem.
\bgn{proof}
If the logarithmic transformation determined by $\Psi$ is trivial, the statement is obvious. 
We assume that the logarithmic transformation is not trivial. 

Let $\omega_T$ be a symplectic form of $\nu T$ which induces the generalized complex structure $\mathcal{J}|_{\nu T}$ and $T$ is symplectic with respect to $\ome_T$.
By Weinstein's neighborhood theorem \cite{Weinstein}, we can take a symplectomorphism: 
$$
\Theta: (\nu T, \omega_T) \rightarrow (D^2\times T^2, \sigma_C), 
$$
where $D^2$ denotes the unit disk $\{\, z_1\in \mathbb{C}\, |\, |z_1|\leq1\,\}$ and $T^2$ is the quotient $\mathbb{C}/\mathbb{Z}^2 \cong \mathbb{R}^2/\mathbb{Z}^2$ with a coordinate $z_2$
and  $\sigma_C= \sqrt{-1}C (dz_1\wedge d\bar{z}_1 + dz_2 \wedge d\bar{z}_2)$ for 
the constant $C = \frac{1}{2}\int_{T}\omega_T$.
Using this identification, the attaching map $\Psi$ can be regarded as a matrix $A_\Psi \in $SL$(3;\mathbb{Z})$. 
Any matrix $P\in $SL$(2;\mathbb{Z})$ induces a self-diffeomorphism $P: T^2\rightarrow T^2$. 
Since the map $\id_{D^2}\times P$ preserves the form $\sigma_C$ and the diffeomorphism type of $M_\Psi$ is determined by the first row of $A_\Psi$, we can assume that $A_\Psi$ is equal to the following matrix: 
$$
\begin{pmatrix}
m & 0 & p \\
0 & 1 & 0 \\
a & 0 & b
\end{pmatrix}
$$
(see \cite{Gompf_Stipsicz}),
where $m,p,a,b$ are integers which satisfy $mb-pa=1$. 
The first row of $A_\Psi$ is $(m,0,p)$ and it follows that $p\neq 0$ since $\Psi$ is not trivial.
Thus we can take $a$ and $b$ satisfying the condition 
\begin{equation}\label{eq:condition} mbs -pa \neq 0\quad \text{\rm for all } s\in[0,1] \end{equation} by replacing 
$(a,b)$ by $(a+ml, b+pl)$ for a suitable integer $l$ if necessary.
Let $D_k $ be the annulus $\{z_1 \in \mathbb{C} \hspace{.3em}|\hspace{.3em} k<|z|\leq 1\}$ for $k \in [0,1]$. 
We define a diffeomorphism $\Psi: D_{\frac{1}{e}} \times T^2 \rightarrow D_{0}\times T^2$ as follows: 
$$
\Psi(r, \theta_1,\theta_2, \theta_3) = (\sqrt{\log(er)}, m\theta_1+a\theta_3, \theta_2, p\theta_1+b\theta_3), 
$$
where $z_1=r \exp(\sqrt{-1}\theta_1)$ and $z_2=\theta_2 + \sqrt{-1}\theta_3$. 
The manifold $M_{\Psi}$ is diffeomorphic to the following manifold: 
$$
M_{\Psi} = (X\bsh\Int(\nu T)) \cup_{\Psi} D^2\times T^2. 
$$
Thus it suffices to construct a twisted generalized complex structure on $M_\Psi$ which satisfies the conditions  (\ref{ali:non-deg,pure spinor on 4-mfd1}), (\ref{ali:non-deg,pure spinor on 4-mfd2}) and (\ref{eq:GCS on 4-mfd}) 

We denote by $\varphi_T:=z_1\phi_T\in \wedge^{\text{even}}D^2\times T^2$ the following form: 
\begin{equation}\label{eq:z1phi local form}
z_1 \exp\left(-\frac{mC}{2}\varrho(|z_1|^2)\frac{dz_1}{z_1}\wedge \frac{d\bar{z}_1}{\bar{z}_1} - bC dz_2\wedge d\bar{z}_2 
+\frac{dz_1}{z_1}\wedge dw_2\right), 
\end{equation}
where $w_2=(C(\frac a2-p)z_2-C(\frac a2+p) \bar{z}_2)$ and 
$\varrho: \mathbb{R}\rightarrow [0,1]$ is a monotonic increasing function which satisfies $\varrho(r)=0$ if $|r|<\frac{1}{2e^2}$ and $\varrho(r)=1$ if $|r|\geq \frac{1}{e^2}$. 
The form $\varphi_T$ satisfies the condition (\ref{ali:non-deg,pure spinor on 4-mfd1}). 
It follows form the condition (\ref{eq:condition})  that the top-degree part 
$\lan\varphi_T \wedge \overline{\varphi_T}\ran$ is not trivial. Thus the condition (\ref{ali:non-deg,pure spinor on 4-mfd2}) holds.
Since ${z_1}^{-1}\varphi_T$ is $d$-closed on $D_0\times T^2$, the form $\varphi_T$ satisfies the condition (\ref{eq:GCS on 4-mfd}).
Thus $\varphi_T$ gives a generalized complex structure on $D^2\times T^2$. 
Denote by $B$ and $\omega$ the real part and the imaginary part of the degree-$2$ part of the form $\phi_T={z_1}^{-1}\varphi_T$, respectively. 
Then it follows from a direct calculation that the pullback $\Psi^\ast \sigma_C$ is equal to $\omega$. 
We take a monotonic decreasing function $\tilde{\varrho}:\mathbb{R}\rightarrow [0,1]$ which satisfies $\tilde{\varrho}(r)=1$ for $|r|<\frac{1}{2}$ and $\tilde{\varrho}(r)=0$ for $|r|\geq 1-\varepsilon$, where $\varepsilon>0$ is a sufficiently small number. 
We define a $2$-form $\tilde{B} \in A^2(M\setminus T)$ by 
\begin{equation*}
\tilde{B} = \begin{cases}
{\Psi^{-1}}^\ast (\tilde{\varrho}(|z_1|^2) B) & \text{on }\nu T\setminus T, \\
0 & \text{on }M\setminus \nu T.  
\end{cases}
\end{equation*}
Then the manifold $M\setminus T$ admits a twisted generalized complex structure $\mathcal{J}^\prime_\Psi$ such that $\exp(\tilde{B}+ \sqrt{-1}\omega_T) \in \Omega_\mathbb{C}^\bullet(\nu T\setminus T)$ is a local section of the canonical bundle $K_{\mathcal{J}^\prime_\Psi}$.  
Since $\varphi_T$ gives a generalized complex structure on $D^2\times T^2$, the form $\exp((\tilde{\varrho}(|z_1|^2)-1)B) \varphi_T$ induces a 
$(-d(\tilde{\varrho}(|z_1|^2)B))$-twisted generalized complex structure. 
Since the $2$-form $\Psi^\ast(\tilde{B}+ \sqrt{-1}\sigma_C)$ is equal to $\tilde{\varrho}(|z_1|^2)B+\sqrt{-1}\omega$, we obtain a twisted generalized complex structure on $M_\Psi$ which satisfies the conditions (\ref{ali:non-deg,pure spinor on 4-mfd1}), (\ref{ali:non-deg,pure spinor on 4-mfd2}) and (\ref{eq:GCS on 4-mfd}).
If $H^3(M_\Psi)=\{0\}$, then there exists a $2$-form 
$\gam\in A^2(M_\Psi)$ such that 
$d\gam =-d(\tilde{\varrho}(|z_1|^2)B)$.
Then $\exp(\gam+(\tilde{\varrho}(|z_1|^2)-1)B) \varphi_T$
This completes the proof of Theorem \ref{th:Goto and Hayano 1}. 
\end{proof}
We can use logarithmic transformations of general multiplicity to obtain generalized complex structures with arbitrary large number of 
connected components of type changing loci.

\bgn{theorem}{\rm \cite{Goto_Hayano}}
Let $(M, \J)$ be a generalized complex $4$-manifold and $T$ a torus with trivial normal bundle. 
We denote by $M'$ the manifold obtained from $M$ by a $C^\infty$ logarithmic transformation on $T$ of multiplicity $0$. 
Then for every $n>\k(\J)$,
$M'$ admits a twisted generalized complex structure $\J'_n$ with $\k(\J'_n) =n$.
In particular, if $H^3(M')=\{0\}$, then $\J_n$ is a generalized complex structure.
\end{theorem}

\bgn{theorem}{\rm \cite{Goto_Hayano}}\label{th:GCS on (2k+1)C P2  l ol C P2}
For every $k,l\geq 0$ and $m\geq 1$, the connected sum 
$(2k+1)\C P^2\#l\ol{\C P^2}$ admits generalized complex structures $\J_m$ with $\k(\J_m)=m$.
\end{theorem}
\bgn{remark} 
Cavalcanti and Gualtieri firstly  constructed a generalized complex structure $\J$ on the connected sum $(2k+1)\C P^2\#l\ol{\C P^2}$
with $\k(\J)=1$ by logarithmic transformations of multiplicity $0$ which does not admit any complex structures and symplectic structures \cite{Cavalcanti_Gualtieri_2006}, \cite{Cavalcanti_Gualtieri_2009}.
\end{remark}
\bgn{remark}
In \cite{Cavalcanti_Gualtieri_2006}
it is pointed out that another  generalized complex surgery is possible, but little detail is
provided. Working out the details of this another surgery one would obtain an alternative proof of Theorem \ref{th:Goto and Hayano 1}
in the special cases where the first row of the matrix $A_\Psi$ is  $(m,0,1).$
\end{remark}

\bgn{remark}
Torres and Yazinski also constructed twisted generalized complex manifolds on several manifolds with arbitrary large 
type changing loci by a different method \cite{Torres_Yazinski}. 
\end{remark}
We shall apply our theorems to these $4$-manifolds obtained by $C^\infty$ logarithmic transformations.
\bgn{theorem}\label{th:unobstruct deformations of log 4-mfds}
Let $(M, \J)$ be a {generalized complex structure}  and $T$ a real symplectic torus with trivial normal bundle. 
We denote by $(M_\Psi,\J)$  a generalized complex $4$-manifold constructed by a $C^\infty$ logarithmic transformation from $(M, \J)$ along the torus $T$
with type changing loci $D$
as in Theorem \ref{th:Goto and Hayano 1}, where we assume that  $H^3(M_\Psi)=\{0\}$.
Then the Lie algebroid cohomology $H^k(\w^\bullet \ol L_{\J_\Psi})$ is isomorphic to  $H^k(M_\Psi\bsh D, \C)$ and
deformations of generalized complex structures of $(M_\Psi,\J_\Psi)$ are unobstructed which are parametrized by an open set of $H^2(M_\Psi\bsh D, \C)$.
\end{theorem}
\bgn{proof}
It follows from (\ref{eq:z1phi local form}) that the degree $2$-part $\ome_\C$ of $\phi_T$ is given by 
\begin{equation}
\left( - bC dz_2\wedge d\bar{z}_2 + C\left( \frac{a}{2} -p \right) \frac{dz_1}{z_1}\wedge dz_2 -C\left( \frac{a}{2}+p \right) \frac{dz_1}{z_1} \wedge d\bar{z}_2\right), 
\end{equation}
if $r=|z_1|<\frac 1{e^2}$. 
On the complement of $D$, the $\ome_\C$ is given by $b+\sqrt{-1}\ome$ for a real $2$-form $b$ and a real symplectic form $\ome$. 
The $2$-forms $\ome_\C$ is written as 
$$
\frac{dz_1}{z_1}\w dw_2- C'dw_2\w d\ol w_2,
$$
where $w_2=C(\frac a2-p)z_2-C(\frac a2+p)\ol z_2,$ and $\dstyle{C'=\frac{bc}{C^2(a^2+4p^2)}}$. 

We define local $C^\infty$ complex coordinates $(w_1, w_2)$ by 
$$
w_1=e^{C' \ol w_2} z_1, \quad w_2=C(\frac a2-p)z_2-C(\frac a2+p)\ol z_2.
$$

Note that $p$ is not equal to $0$.
Then we have 
$$
\ome_\C=\frac{dw_1}{w_1}\w dw_2.
$$
Thus $\ome_\C$ is a $C^\infty$ logarithmic symplectic structure on a neighborhood of $D$.
 $\J_{\Psi}$ is induced from $\dstyle{\phi_T:=e^{\kappa+\ome_\C}}$, where 
 $\kappa:=\gam+(\tilde{\varrho}(|z_1|^2)-1)B$ is a $d$-closed real $2$-form on $M_\Psi$. 
From Remark \ref{re:Lie algebroid iso b-fields}, the Lie algebroid cohomology groups are invariant under the action of $d$-closed $b$-fields.
 The action of $d$-closed $b$-fields also preserves the unobstructedness of deformations of {generalized complex structure}.
Thus the results follow from 
 Theorem \ref{th:log.cohomology} and Theorem \ref{th:unobst. log GCS}.
\end{proof}
\bgn{theorem}\label{unobst. def of (2k-1)CP2 (10k-1)ol CP2}
Let $\J_m$ be the {generalized complex structure} on the connected sum $M:=(2k-1)\C P^2\#(10k-1)\ol{\C P^2}$ as in Theorem \ref{th:GCS on (2k+1)C P2  l ol C P2} and $D$ the type changing loci of $\J_m$. Then the Lie algebroid cohomology $H^k(\w^\bullet \ol L_\J)$ is isomorphic to  $H^k(M\bsh D, \C)$ and
deformations of the generalized complex structure $\J_m$ are unobstructed which are parametrized by an open set of $H^2(M\bsh D, \C)$.
\end{theorem}
\bgn{proof}
Let $E(1)$ be the blown up $\C P^2$ at $9$ points of intersection of two generic cubic hypersurfaces in $\C P^2$. 
The manifold $E(1)$ is an elliptic fibration over $\C P^1$. 
We denote by $E(k)$ the fibre sum of $k$ copies of $E(1)$. 
Then a $C^\infty$ logarithmic transformation of multiplicity $0$ along a smooth fibre of $E(k)$
gives a manifold $M:=(2k-1)\C P^2\#(10k-1)\ol{\C P^2}$.
As in \cite{Goto_Hayano}, we can apply 
a $C^\infty$ logarithmic transformations of multiplicity $1$ on fibres of fishtail
neighborhood at $m$ times so that the diffeomorphism type of 
$M:=(2k-1)\C P^2\#(10k-1)\ol{\C P^2}$ is not changed. 
The procedure gives the {generalized complex structure} $\J_m$ on $M:=(2k-1)\C P^2\#(10k-1)\ol{\C P^2}$ such that $\kappa(\J_m)=m$. 
Thus it follows from Theorem \ref{th:unobstruct deformations of log 4-mfds} that $\J_m$ is induced from a logarithmic symplectic structure acted by a $d$-closed $b$-field.
The result follows from Theorem \ref{th:log.cohomology} and Theorem \ref{th:unobst. log GCS}. 
\end{proof}
\bgn{remark}
Let $D=\cup_{i=1}^m T_i$ be the type changing loci of 
{generalized complex structure} $\J_m$ in Theorem \ref{unobst. def of (2k-1)CP2 (10k-1)ol CP2}. 
Then $E(n)\bsh D\cong M\bsh D$ and we obtain an exact sequence
$$
H^2_D((E(k))\overset i\to H^2(E(n))\to H^2(E(n)\bsh D)\to H^3_D(E(n))\to 0,
$$
where $H^i_D(E(k))$ denotes the local cohomology of $E(k)$ with support $D$. 
By the duality,  we have $H^i_D(E(k))\cong H^{i-2}(D)$. 
The image of the map $i:H^2_D((E(k))\to H^2(E(n))$ is a one dimensional subspace of $H^2(E(k))$ 
which is generated by the $1$-st Chern class of the line bundle given by a fibre of $E(n)$. 
Thus we have $\dim H^2(M\bsh D) =\dim H^2(E(k))+\dim H^1(D)-1.$
Since $\dim H^2(E(k))=12k-2$, it follows that $\dim H^2(M\bsh D) =12k+2m-3$.
\end{remark}



\end{document}